\documentclass{amsart}
\usepackage{graphicx}
\usepackage{amsmath}
\usepackage{amsthm}
\usepackage{amsfonts}
\usepackage{amssymb,amscd}

\usepackage{color}

\usepackage[all]{xy}

\newtheorem{theorem}{Theorem}[section]

\newtheorem{lem}[theorem]{Lemma}
\newtheorem{proposition}[theorem]{Proposition}
\newtheorem{example}[theorem]{Example}
\theoremstyle{definition}
\newtheorem{definition}[theorem]{Definition}
\theoremstyle{remark}
\newtheorem{remark}[theorem]{Remark}
\numberwithin{equation}{section}

\newcommand{\K}{\mathbb K}

\newcommand{\A}{\mathcal{A}}

\begin{document}

\title[ Generalization of $n$-ary Nambu algebras and beyond]
{ Generalization of $n$-ary Nambu algebras and beyond}%
\author{H. Ataguema, A. Makhlouf and S. Silvestrov}%
\address{Abdenacer Makhlouf and Hammimi Ataguema, Universit\'{e} de Haute Alsace,  Laboratoire de Math\'{e}matiques, Informatique et Applications,
4, rue des Fr\`{e}res Lumi\`{e}re F-68093 Mulhouse, France}%
\email{Abdenacer.Makhlouf@uha.fr}
\address{Sergei Silvestrov, Centre for Mathematical Sciences,  Lund University, Box
   118, SE-221 00 Lund, Sweden}
\email{sergei.silvestrov@math.lth.se}

\thanks {
This work was
partially supported by the Crafoord Foundation, The
Swedish Foundation for International Cooperation in Research and Higher Education (STINT), The Swedish research links program of SIDA foundation and the Swedish Research Council, The Royal Swedish Academy of Sciences,
The Royal Physiographic Society in Lund, The
European network Liegrits,  Mulhouse University and Lund University.}

\date{November 2008}
%
\begin{abstract}
The aim of this paper is to introduce $n$-ary Hom-algebra structures
generalizing the $n$-ary algebras of Lie type enclosing   $n$-ary
Nambu algebras, $n$-ary Nambu-Lie algebras, $n$-ary Lie algebras, and
$n$-ary algebras of associative type enclosing $n$-ary totally
associative and $n$-ary partially associative algebras. Also, we
provide a way to construct examples starting from an $n$-ary algebra
and an $n$-ary algebras endomorphism. Several examples could be
derived using this process.

\end{abstract}
\maketitle

\section*{Introduction}
In this paper, we introduce generalizations of $n$-ary algebras of
Lie type and associative type by twisting the identities using
linear maps. The $n$-ary algebraic structures and in particular
ternary algebraic structures appeared more or less naturally in
various domains of theoretical and mathematical physics and data
processing. Indeed, theoretical physics progress of quantum
mechanics and the discovery of the Nambu mechanics  (1973) (see \cite{Nambu}), as well
as a work of S. Okubo (see \cite{Okubo}) on Yang-Baxter equation gave  impulse to a
significant development on $n$-ary algebras.  The $n$-ary operations
appeared first through  cubic matrices which were introduced in the
nineteenth century by Cayley. The cubic matrices were considered
again and generalized by Kapranov, Gelfand, Zelevinskii in 1994 (see
\cite{KapranovGelfandZelinski}) and Sokolov in 1972 (see
\cite{Sokolov}).
Another recent motivation to study $n$-ary operation comes from string
theory and M-branes involving naturally an algebra with
ternary operation called Bagger-Lambert algebra \cite{BL2007}.
Hundred of papers are dedicated to Bagger-Lambert algebra.  For
other physical applications (see \cite{Kerner, Kerner2, Kerner3,Kerner4}).

The first conceptual generalization of binary algebras was the
ternary algebras introduced by Jacobson (see \cite{Jacobson}).
In connection with problems from Jordan theory and quantum mechanics,
he defined the Lie triple systems. A Lie triple system consists of a
space of linear operators on vector space $V$ that is closed under
the ternary bracket $[x,y,z]_T=[[x,y],z]$, where $[x,y]=x y-y x$.
Equivalently, the Lie triple system may be viewed as a subspace of
the Lie algebra closed relative to the ternary product. A Lie triple
system arose also in the study of symmetric spaces (see \cite{loosSym}).
More generally, we distinguish two kinds of generalizations of binary
Lie algebras. Firstly, $n$-ary Lie algebras in which the Jacobi identity is
generalized by considering a cyclic summation over
$\mathcal{S}_{2n-1}$ instead of $\mathcal{S}_3$, (see \cite{Hanlon}
\cite{Michor}) and secondly $n$-ary Nambu algebras in which the fundamental
identity generalizes the fact that the adjoint maps are derivations.
The fundamental identity appeared first in Nambu mechanics
\cite{Nambu}, the abstract definition of $n$-ary Nambu algebras or
$n$-ary Nambu-Lie algebras (when the bracket is skew symmetric) was
given by Fillipov  in 1985 (see \cite{Filippov}). See also
(\cite{Takhtajan,Takhtajan1}) for the algebraic formulation of the
Nambu mechanics. The Leibniz $n$-ary algebras were introduced and
studied in \cite{CassasLodayPirashvili}. While the basic concepts of
ternary Hopf algebra were introduced in \cite{Santana} and \cite{Duplij}, see also
\cite{GozeRaush}.

On the other hand, $n$-ary algebras of associative type were studied
by Lister, Loos, Myung and Carlsson (see \cite{Carlsson,Lister,Loos,MyungTernAlg73}).
The $n$-ary operations of associative type lead to two principal classes of "associative" $n$-ary algebras,  totally associative $n$-ary algebras and partially associative $n$-ary  algebras.
Also they admit  some variants. Typical examples of totally associative $n$-ary algebras for $n=3$ are given by
subspaces of an associative algebra that are closed relative to the ternary product $(xyz)\mapsto xyz$.
This natural ternary product is related to the important ternary operation introduced by Hestenes (see \cite{Hestenes}) defined on the linear space of rectangular  matrices $A,B,C\in \mathcal{M}_{m,n}$ with complex entries by $AB^{\ast }C$ where $B^{\ast}$ is the conjugate transpose matrix of $B$. This operation is strictly speaking not a ternary algebra product on $\mathcal{M}_{m,n}$ as it is linear on the first and the third arguments but conjugate-linear on the second argument. It satisfies identities, sometimes referred to as identities of total associativity of second kind, only slightly different from the identities of totally
associative algebras. The totally associative ternary algebras are also sometimes called associative triple systems. The cohomology of totally associative $n$-ary algebras was studied by
Carlsson through the embedding (see \cite{Carlsson1}). In
\cite{AM2007}, the first and second authors  extended the 1-parameter
formal deformation theory to ternary algebras of associative type
and in \cite{AM2008} they discussed their cohomologies adapted to deformations, see also \cite{Hoffbeck,GozeRemm,Remm}.

The Hom-algebras structures arose first in quasi-deformation of Lie
algebras of vector fields. These quasi-deformations lead to quasi-Lie algebras, a
generalized Lie algebra structure in which the skew-symmetry and Jacobi conditions are twisted. The first examples were concerned with $q$-deformations of the Witt and Virasoro algebras  (see \cite{AizawaSaito,ChaiElinPop,ChaiKuLukPopPresn,ChaiIsKuLuk,ChaiPopPres,CurtrZachos1,
DaskaloyannisGendefVir,Kassel1,
LiuKeQin,Hu}). Motivated by these and the new examples arising as
application of the general quasi-deformation construction of
\cite{HLS,LS1,LS2} on the one hand, and the desire to be able to
treat within the same framework such well-known generalizations of
Lie algebras as the color and super Lie algebras on the other hand,
quasi-Lie algebras and subclasses of quasi-Hom-Lie algebras and
Hom-Lie algebras were introduced in \cite{HLS,LS1,LS2,LS3}. In the
subclass of Hom-Lie algebras skew-symmetry is untwisted, whereas the
Jacobi identity is twisted by a single linear map and contains three
terms as in Lie algebras, reducing to ordinary Lie algebras when the
twisting linear map is the identity map.

The notion of Hom-associative algebras generalizing
associative algebras to a situation where
associativity law is twisted by a linear map was introduced by the second and the third authors in \cite{MS}, where it was shown in particular that the
commutator product, defined using the
multiplication in a Hom-associative algebra, leads
naturally to a Hom-Lie algebra. We introduced also
Hom-Lie-admissible algebras and more general
$G$-Hom-associative algebras with subclasses of
Hom-Vinberg and pre-Hom-Lie algebras,
generalizing to the twisted situation
Lie-admissible algebras, $G$-associative
algebras, Vinberg and pre-Lie algebras
respectively, and shown that for these classes of
algebras the operation of taking commutator leads
to Hom-Lie algebras as well.
The enveloping algebras of Hom-Lie algebras were discussed in \cite{Yau:EnvLieAlg}.
It is important challenging problem to develop further the proper notion and theory of
universal enveloping algebras for general Hom-algebras and for general
quasi-Hom-Lie and Quasi-Lie algebras.
The fundamentals of the formal
deformation theory and associated cohomology
structures for Hom-Lie algebras have been
considered recently by the second and the third authors in
\cite{HomDeform}. Simultaneously, D. Yau has
developed elements of homology for Hom-Lie
algebras in \cite{Yau:HomolHomLie}. In \cite{HomHopf} and \cite{HomAlgHomCoalg}, we
developed the theory of Hom-coalgebras and related structures.
We introduced Hom-coalgebra
structure, leading to the notions of Hom-bialgebra
and Hom-Hopf algebra, proved some fundamental properties and provided examples.
Also, we defined the concept of Hom-Lie
admissible Hom-coalgebra generalizing the admissible coalgebra introduced in \cite{GR},
and provide their classification based on subgroups of the symmetric
group.

This paper is organized as follows.
Section \ref{sec1:naryHomalgLietype} is dedicated to a generalization
 of the $n$-ary algebras of Lie type to Hom structures. We introduce
the notions of $n$-ary Hom-Nambu algebra, $n$-ary Hom-Nambu-Lie
algebra and $n$-ary Hom-Lie algebra. Also, we provide a process to
construct them starting from $n$-ary Nambu algebra (resp. $n$-ary
Nambu-Lie algebra or $n$-ary Lie algebra) and an algebra
endomorphism of this $n$-ary algebra. In Section~\ref{sec2:naryHomAlgAssType}, we
introduce the $n$-ary Hom-algebras of associative type, and provide a
procedure for construction of such $n$-ary Hom-algebras from $n$-ary algebras of the same
associative type and their $n$-ary algebra endomorphisms. We consider
the totally Hom-associative $n$-ary algebras and partially
Hom-associative $n$-ary algebras. Their variants, a weak totally
associative $n$-ary algebras and the alternate partially associative
$n$-ary algebras, could be also considered in the same scheme. We
provide several examples, using this procedure. In Section \ref{sec3:tensorproducts}, we
present some results concerning the tensor product of $n$-ary Hom-algebras.

\section{The $n$-ary Hom-algebras of Lie Type} \label{sec1:naryHomalgLietype}
Throughout this paper, we will for simplicity of exposition assume that
$\K$ is an algebraically closed field of
characteristic zero, even though for most of the general definitions and results in the paper
this assumption is not essential.
Let $V$ be a $\K$-vector space.

\subsection{Definitions}
In this section, we introduce the definition of $n$-ary   Hom-Nambu algebras,
 $n$-ary   Hom-Nambu-Lie algebras and $n$-ary   Hom-Lie algebras,
 generalizing the usual  $n$-ary   Nambu algebras,
 $n$-ary   Nambu-Lie algebras (called also Filippov $n$-ary algebras)
 and $n$-ary   Lie algebras.
\begin{definition}
An $n$-\emph{ary  Hom-Nambu algebra} is a triple
$(V,[\cdot,\dots,\cdot],\alpha)$, consisting of a vector space $V$,
an $n$-linear map $[\cdot,\cdots,\cdot]:V^{\times n}\rightarrow V$ and a
family $\alpha=(\alpha_{i})_{i=1,\cdots,n-1}$ of linear maps
$\alpha_{i}:V\rightarrow V$, $i=1,\cdots ,n-1$ satisfying
\begin{eqnarray}\label{NambuIdentity}
[ \alpha_1(x_{1}),..., \alpha_{n-1}(x_{ n-1}),
[x_{n},...,x_{2n-1}]]=\\
\sum_{i=n}^{2n-1}{[\alpha_1(x_{n}),...,\alpha_{i-n}(x_{i-1}), [
x_{1},\cdots,x_{n-1},x_{i}], \alpha_{i-n+1}(x_{i+1})
...,\alpha_{n-1}(x_{2n-1})]}\nonumber
\end{eqnarray}%
for all $(x_1,\cdots,x_{2n-1}) \in V^{2n-1}.$
\end{definition}
We call the condition (\ref{NambuIdentity}) the $n$-ary Hom-Nambu
identity.

In particular, the  ternary Hom-Nambu algebras are defined by the
following Hom-Nambu identity
\begin{eqnarray}\label{TernaryNambuIdentity}
[ \alpha_1(x_{1}),\alpha_{2}(x_{ 2}), [x_{3},x_{4},x_{5}]]=  \quad \quad \quad\quad \quad \quad\quad \quad \quad \\
\quad \quad \quad [ [x_1,x_2,x_3 ],\alpha_1(x_{4}),\alpha_{2}(x_{5})
] + [ \alpha_{1} (x_3),[ x_1,x_2,x_{4}] , \alpha_{2}(x_{5})]+
[\alpha_1(x_3),\alpha_{2}(x_{4}),[x_1,x_2,x_{5}]]. \nonumber
\end{eqnarray}%

\begin{remark}
Let $(V,[\cdot,\dots,\cdot],\alpha)$ be an $n$-ary Hom-Nambu algebra
where $\alpha=(\alpha_{i})_{i=1,\cdots,n-1}$. For any
$x=(x_1,\ldots,x_{n-1})\in V^{n-1}$, $\alpha
(x)=(\alpha_1(x_1),\ldots,\alpha_{n-1}(x_{n-1}))\in V^{n-1}$ and
$y\in V$, let $L_x$ be a linear map on $V$,
defined by  $$L_{x}(y)=[x_{1},\cdots,x_{n-1},y].$$

Then the Hom-Nambu identity may be written as
\begin{equation*}
L_{\alpha (x)}( [x_{n},\dots,x_{2n-1}])=
\sum_{i=n}^{2n-1}{[\alpha_1(x_{n}),\dots,\alpha_{i-n}(x_{i-1}),
L_{x}(x_{i}), \alpha_{i-n+1}(x_{i+1}),\dots,\alpha_{n-1}(x_{2n-1})].}
\end{equation*}%
\end{remark}
\begin{remark}
When the maps $(\alpha_{i})_{i=1,\cdots,n-1}$ are all identity maps,
one recovers the classical $n$-ary   Nambu algebras.
In the special case of classical $n$-ary Nambu algebra, it's known
also as the fundamental identity or Filippov identity.
\end{remark}

\begin{definition}
An $n$-ary  Hom-Nambu algebra
$(V,[\cdot,\dots,\cdot],\alpha)$ where $\alpha=(\alpha_{i})_{i=1,\cdots,n-1}$ is called
$n$-\emph{ary  Hom-Nambu-Lie algebra} if the bracket is
skew-symmetric that is
\begin{equation}
[x_{\sigma (1)}, \cdots, x_{\sigma (n)}]=Sgn(\sigma )[x_{1},\cdots, x_{n}],\quad \forall \ \sigma \in \mathcal{S}_{n}\text{
and }\forall \ x_{1},\cdots,x_{n}\in V
\end{equation}%
where $\mathcal{S}_{n}$ stands for the permutation group on
$m$ elements.
\end{definition}

\begin{definition}
An $n$-\emph{ary  Hom-Lie algebra} is  a triple
$(V,[\cdot,\dots,\cdot],\alpha)$, consisting of a vector space
$V$, a skew-symmetric map $[\cdot,\cdots,\cdot]:V^{\times n}\rightarrow V$
and a family $\alpha=(\alpha_{i})_{i=1,\cdots,n-1}$ of linear maps $\alpha_{i}:V\rightarrow V$, $i=1,\cdots
,n-1$ satisfying $n$-\emph{ary Hom-Jacobi identity}
\begin{equation}\label{nLieIdentity}
\sum_{\sigma \in \mathit{S}_{2n-1}}{ Sgn(}\sigma
)[\alpha_1(x_{\sigma (1)}),..., \alpha_{n-1}(x_{\sigma
(n-1)}),[x_{\sigma (n)},...,x_{\sigma (2n-1)}]]=0
\end{equation}%
for all $(x_1,\cdots,x_{2n-1}) \in V^{2n-1}$.
\end{definition}
In particular, the ternary  Hom-Lie algebras are defined by the
following ternary Hom-Jacobi identity
\begin{equation*}
{\sum_{\sigma \in \mathit{S}_{5}}{ Sgn(}\sigma )[\alpha_1(x_{\sigma
(1)}),\alpha_{2}(x_{\sigma (2)}),[x_{\sigma (3)},x_{\sigma
(4)},x_{\sigma (5)}]]=0.}
\end{equation*}%

The morphisms of $n$-ary Hom-algebras of Lie type are defined as
follows.

\begin{definition}
Let $(V,[\cdot,\dots,\cdot],\alpha)$ and
$(V',[\cdot,\dots,\cdot]',\alpha')$ be two $n$-ary Hom-Nambu
algebras (resp. $n$-ary Hom-Nambu-Lie algebras, $n$-ary Hom-Lie
algebras) where $\alpha=(\alpha_{i})_{i=1,\cdots,n-1}$ and
$\alpha'=(\alpha'_{i})_{i=1,\cdots,n-1}$. A linear map $\rho:
V\rightarrow V'$ is an  $n$-ary Hom-Nambu algebras morphism (resp.
$n$-ary Hom-Nambu-Lie algebras morphism, $n$-ary Hom-Lie algebras
morphism) if it satisfies
\begin{eqnarray*}\rho ([x_{1},\cdots,x_{n}])&=&
[\rho (x_{1}),\cdots,\rho (x_{n})]'\\
\rho \circ \alpha_i&=&\alpha'_i\circ \rho \quad \forall \ i=1,\dots,n-1.
\end{eqnarray*}
\end{definition}

\subsection{Construction of $n$-ary   Hom-algebras of Lie type}
The following theorem provides a way to construct an $n$-ary
Hom-Nambu algebra (resp. $n$-ary Hom-Nambu-Lie algebra, $n$-ary
Hom-Lie algebra) starting from an  $n$-ary Nambu algebra (resp. $n$-ary
Nambu-Lie algebra, $n$-ary Lie algebra) and an $n$-ary algebras
endomorphism.
\begin{theorem} \label{thm:naryNambumorphism}
Let $(V,m)$ be an $n$-ary Nambu algebra (resp. $n$-ary Nambu-Lie
algebra,  $n$-ary Lie algebra)  and let $\rho : V\rightarrow V$ be
an $n$-ary  Nambu (resp.  $n$-ary Nambu-Lie, $n$-ary Lie) algebras
endomorphism.

We set $m_\rho=\rho\circ m$ and
$\widetilde{\rho}=(\rho,\cdots,\rho)$. Then
$(V,m_\rho,\widetilde{\rho})$ is an $n$-ary Hom-Nambu
(resp. $n$-ary Hom-Nambu-Lie, $n$-ary Hom-Lie) algebra.

Moreover, suppose that  $(V',m')$ is another $n$-ary Nambu algebra
(resp.  $n$-ary Nambu-Lie algebra, $n$-ary Lie algebra)
and $\rho ' : V'\rightarrow V'$ is a $n$-ary Nambu
(resp.  $n$-ary Nambu-Lie, $n$-ary Lie) algebra endomorphism. If
$f:V\rightarrow V'$ is a $n$-ary Nambu algebra morphism (resp.
$n$-ary Nambu-Lie algebra morphism, $n$-ary Lie algebra morphism)
that satisfies $f\circ\rho=\rho'\circ f$ then
$$f:(V,m_\rho,\widetilde{\rho})\longrightarrow (V',m'_{\rho '},\widetilde{\rho '})
$$
is an  $n$-ary Hom-Nambu algebras morphism (resp.  $n$-ary
Hom-Nambu-Lie algebras morphism, $n$-ary Hom-Lie algebras morphism).
\end{theorem}

\begin{proof}

We show that $(V,m_\rho,\widetilde{\rho})$ satisfies the Hom-Nambu identity.
Indeed
\begin{align*}
m_\rho(\rho(x_1),\cdots,\rho(x_{n-1}),m_\rho (x_{n},\cdots x_{2n-1}))&=
\rho\circ m(\rho(x_1),\cdots,\rho(x_{n-1}),\rho\circ m (x_{n},\cdots x_{2n-1}))\\
\ & =\rho^2\circ m(x_1,\cdots,x_{n-1}, m (x_{n},\cdots x_{2n-1})).
\end{align*}
On the other hand,
\begin{eqnarray*}
\sum_{i=1}^{n}{m_\rho(\rho(x_{n}),...,\rho(x_{i-1}), m_\rho(
x_{1},\cdots,x_{n-1},x_{i}), \rho(x_{i+1})
...,\rho(x_{2n-1}))}\\ =
\sum_{i=1}^{n}{\rho \circ m(\rho(x_{n}),...,\rho(x_{i-1}), \rho \circ m(
x_{1},\cdots,x_{n-1},x_{i}), \rho(x_{i+1})
...,\rho(x_{2n-1}))}\\
=\rho^2(\sum_{i=1}^{n}{  m(x_{n},...,x_{i-1},  m(
x_{1},\cdots,x_{n-1}),x_{i}, x_{i+1}
...,x_{2n-1}))}.
\end{eqnarray*}
Therefore the Nambu identity implies the Hom-Nambu identity for this
$n$-ary algebra.
The skew-symmetry and the Hom-Jacobi identity
are proved similarly.
The second assertion follows from
$$ f\circ m_\rho = f\circ \rho \circ m
= \rho ' \circ f \circ m = \rho ' \circ m' \circ f
=  m'_{\rho '} \circ f. $$
\end{proof}

\begin{example}
 An algebra $V$ consisting of polynomials or possibly of other differentiable functions in $3$ variables $x_{1},x_{2},x_{3},$ equipped with well-defined bracket multiplication given by the functional jacobian
 $J(f)=(\frac{\partial f_i}{\partial x_j})_{1\leq i,j \leq 3}$:
\begin{equation} \label{JacternaryNambuLie}
\lbrack f_{1},f_{2},f_{3}\rbrack=\det \left(
\begin{array}{ccc}
\frac{\partial f_{1}}{\partial x_{1}} & \frac{\partial f_{1}}{\partial
x_{2}} &
\frac{\partial f_{1}}{\partial x_{3}} \\
\frac{\partial f_{2}}{\partial x_{1}} & \frac{\partial f_{2}}{\partial
x_{2}} &
\frac{\partial f_{2}}{\partial x_{3}} \\
\frac{\partial f_{3}}{\partial x_{1}} & \frac{\partial f_{3}}{\partial
x_{2}} &
\frac{\partial f_{3}}{\partial x_{3}}%
\end{array}%
\right),
\end{equation}%
is a ternary  Nambu-Lie algebra. By considering a ternary Nambu-Lie algebra endomorphism of such algebra, we construct a ternary  Hom-Nambu-Lie algebra.
Let $\gamma(x_1,x_2,x_3)$ be a polynomial or more general differentiable transformation of three variables mapping elements of $V$ to elements of $V$ and such that $\det J (\gamma) = 1$. Let $\rho_\gamma: V \mapsto V$ be the composition transformation defined by
$f\mapsto f\circ \gamma $ for any $f \in V$.
By the general chain rule for composition of transformations of several variables,
\begin{eqnarray*}
J(\rho_\gamma (f)) = J(f\circ \gamma)=(J(f)\circ \gamma) J(\gamma)=\rho_\gamma (J(f))J(\gamma), \\
\det J(\rho_\gamma (f)) = \det (J(f)\circ \gamma) \det J(\gamma)= \det \rho_\gamma (J(f))\det J(\gamma).
\end{eqnarray*}
Hence, for any transformation $\gamma$ with $\det J(\gamma)=1$, the composition transformation
$\rho_\gamma$ defines an endomorphism of the ternary  Nambu-Lie algebra with ternary product  \eqref{JacternaryNambuLie}.
Therefore, by Theorem \ref{thm:naryNambumorphism}, for any such transformation $\gamma$, the triple
$\left(V,\lbrack \cdot,\cdot,\cdot \rbrack_\gamma=\rho_\gamma\circ \lbrack \cdot,\cdot,\cdot \rbrack,\widetilde{\rho_\gamma}\right)$
is a ternary Hom-Nambu-Lie algebra.
\end{example}

\section{The $n$-ary Hom-Algebras of Associative Type} \label{sec2:naryHomAlgAssType}
In this Section, we introduce the $n$-ary  totally Hom-associative
algebras and $n$-ary    partially Hom-associative algebras
generalizing the classical $n$-ary  totally associative algebras and
$n$-ary    partially associative algebras.
 We also
provide a process to construct them starting from an  $n$-ary
totally associative algebra or $n$-ary    partially associative
algebra and an algebra endomorphism of such $n$-ary algebra.

\subsection{Definitions}
We introduce the definitions of $n$-ary  totally Hom-associative
algebras and $n$-ary    partially Hom-associative algebras.
\begin{definition}
An $n$-ary  totally Hom-associative algebra is a triple
$(V,m,\alpha)$, consisting of a vector space $V$, a $n$-linear map
$m:V^{\times n}\rightarrow V$ and a family
$\alpha=(\alpha_{i})_{i=1,\cdots,n-1}$ of linear maps
$\alpha_{i}:V\rightarrow V$, $i=1,\cdots ,n-1$ satisfying
\begin{equation}\label{TotAssIdentity}
\begin{array}{c}
m(\alpha_1(x_{1}),\cdots,\alpha_{n-1}(x_{n-1}),m(x_{n},\cdots,x_{2n-1}))=\\
\vdots
\\
m(\alpha_1(x_{1}),\cdots,\alpha_{i}(x_{i}),m(x_{i+1},\cdots,x_{i+n}),
\alpha_{i+1}(x_{i+1+n}),\cdots,\alpha_{n-1}(x_{2n-1}))=\\
\vdots \\
 m(m(x_{1},\cdots,x_{n}),\alpha_1(x_{n+1}),\cdots,\alpha_{n-1}(x_{2n-1})).
\end{array}%
\end{equation}%
\end{definition}

\begin{definition}
  An \emph{$n$-ary  partially Hom-associative algebra} is a triple
$(V,m,\alpha)$, consisting of a vector space $V$, a $n$-linear map
$m:V^{\times n}\rightarrow V$ and a family
$\alpha=(\alpha_{i})_{i=1,\cdots,n-1}$ of linear maps
$\alpha_{i}:V\rightarrow V$, $i=1,\cdots ,n-1$ satisfying
\begin{equation}\label{ParAssIdentity}
\sum_{i=0}^{n-1}{
m(\alpha_1(x_{1}),\cdots,\alpha_{i}(x_{i}),m(x_{i+1},\cdots,x_{i+n}),}
\alpha_{i+1}(x_{i+1+n}),\cdots,\alpha_{n-1}(x_{2n-1}))=0.
\end{equation}
\end{definition}
\begin{remark}
When the maps $(\alpha_{i})_{i=1,\cdots,n-1}$ are all identity maps
then one recovers the classical $n$-ary  totally or partially
associative algebras.
\end{remark}
\begin{remark}
These generalizations to $n$-ary Hom-algebra structures could be extend obviously to the variants of
$n$-ary algebras of associative type. An $n$-ary weak totally
Hom-associative algebras, is given by the identity
\begin{equation*}
m(\alpha_1(x_{1}),\cdots,\alpha_{n-1}(x_{n-1}),m(x_{n},\cdots,x_{2n-1}))=
 m(m(x_{1},\cdots,x_{n}),\alpha_1(x_{n+1}),\cdots,\alpha_{n-1}(x_{2n-1})).
\end{equation*}%
Also, the $n$-ary alternate partially Hom-associative algebras,
where some signs in the identity (\ref{ParAssIdentity}) are moved to
minus, could be defined in the same way.
\end{remark}

The morphisms of $n$-ary Hom-algebras of associative type are
defined as follows.
\begin{definition}
Let $(V,m,\alpha)$ and $(V',m',\alpha')$ be two $n$-ary totally
Hom-associative algebras (resp. $n$-ary partially Hom-associative
algebras) where $\alpha=(\alpha_{i})_{i=1,\cdots,n-1}$ and
$\alpha'=(\alpha'_{i})_{i=1,\cdots,n-1}$. A linear map $\rho:
V\rightarrow V'$ is an  $n$-ary totally Hom-associative algebras
morphism (resp. partially Hom-associative algebras morphism) if it
satisfies
\begin{eqnarray*}\rho (m(x_{1},\cdots,x_{n}))=
m'(\rho (x_{1}),\cdots,\rho (x_{n}))\quad \text{and} \quad  \rho
\circ \alpha_i=\alpha'_i\circ \rho \quad \forall \ i=1,\dots, n-1.
\end{eqnarray*}
\end{definition}

\subsection{Construction of $n$-ary   Hom-algebras of Associative type}
The following theorem provide a way to construct an $n$-ary totally
Hom-associative algebra (resp. $n$-ary partially Hom-associative
algebra) starting from an $n$-ary totally associative (resp. $n$-ary
partially associative algebra) and $n$-ary algebras endomorphism.
\begin{theorem}\label{thm2}
Let $(V,m)$ be an $n$-ary totally associative algebra (resp. $n$-ary
partially associative algebra)  and let $\rho : V\rightarrow V$ be
an $n$-ary  totally Hom-associative (resp. $n$-ary partially
Hom-associative) algebra endomorphism. We set $m_\rho=\rho\circ m$ and
$\widetilde{\rho}=(\rho,\cdots,\rho)$. Then
$(V,m_\rho,\widetilde{\rho})$
  is an $n$-ary totally
Hom-associative algebra (resp. $n$-ary partially Hom-associative
algebra).

 Moreover, suppose that  $(V',m')$ is another $n$-ary totally
associative algebra (resp. $n$-ary partially associative algebra)
and $\rho ' : V'\rightarrow V'$ is a $n$-ary totally associative
 (resp. $n$-ary partially associative) algebras endomorphism. If
$f:V\rightarrow V'$ is a $n$-ary totally associative algebra
morphism (resp. $n$-ary partially associative algebra morphism) that
satisfies $f\circ\rho=\rho'\circ f$,  then
$$f:(V,m_\rho,\widetilde{\rho})\longrightarrow (V',m'_{\rho '},\widetilde{\rho '})
$$
is an $n$-ary totally Hom-associative algebras morphism (resp.
$n$-ary partially Hom-associative algebras morphism).
\end{theorem}

\begin{proof}
We show that $(V,m_\rho,\widetilde{\rho})$ satisfies the identity
(\ref{TotAssIdentity}). Indeed, for $i=1,\cdots,n-1$ we have
\begin{align*}
m_\rho(\rho(x_{1}),\cdots,\rho(x_{i}),m_\rho(x_{i+1},\cdots,x_{i+n}),
\rho(x_{i+1+n}),\cdots,\rho(x_{2n-1}))=\\
\rho\circ m(\rho(x_{1}),\cdots,\rho(x_{i}),\rho\circ
m(x_{i+1},\cdots,x_{i+n}), \rho(x_{i+1+n}),\cdots,\rho(x_{2n-1}))=\\
\rho^2( m(x_{1},\cdots,x_{i}, m(x_{i+1},\cdots,x_{i+n}),
x_{i+1+n},\cdots,x_{2n-1})).
\end{align*}
Therefore the $n$-ary total associativity  identity implies the $n$-ary
total Hom-associativity  identity for this $n$-ary algebra.
The proof for  $n$-ary partially Hom-associative algebras is similar
The second assertion is proved in the same way as the second part of Theorem \ref{thm:naryNambumorphism}.
\end{proof}

\subsection{Examples}
\begin{example}
In dimension 2, there exists a partially associative ternary algebra defined with respect to a basis
$\{ e_1,e_2\}$ as follows
\begin{equation*}
m(e_{1}, e_{1}, e_{1})=e_{2}, \quad m(e_{i}, e_{j}, e_{k})=0, \text{ when } 2\in \{ i,j,k\}.
\end{equation*}
Its automorphisms are given with respect to the same basis by the matrices of the form :
\noindent$\left(
\begin{array}{cc}
a & 0 \\
b & a^3%
\end{array}
\right)$ with $a,b\in \K$ and $a\neq 0$.

We deduce using Theorem \ref{thm2} the 2-dimensional partially
Hom-associative ternary algebra $(\K^2,\widetilde{m},\rho )$ given
with respect to a basis  $\{ e_1,e_2\}$, by the following nontrivial
product
\begin{equation*}
\widetilde{m}(e_{1}, e_{1}, e_{1})=a^3\ e_{2},
\end{equation*}
and the linear map
\begin{eqnarray*}
\rho (e_1)&=&a\ e_1 +b\ e_2\\
\rho (e_2)&=&a^3\ e_2
\end{eqnarray*}
\end{example}

\begin{example}
We consider now  the 2-dimensional totally associative ternary
algebra $(\K^2,m)$ given with respect to a basis $e_1,e_2$ of
$\K^2$, by the following products :

\begin{equation*}
\begin{array}{ll}
\begin{array}{lll}
m(e_{1}, e_{1}, e_{1})& = & e_{1}\\
m(e_{1}, e_{1}, e_{2})& = & e_{2} \\
m(e_{1}, e_{2}, e_{2})& = & e_{1}+e_{2} \\
m(e_{2}, e_{1}, e_{1})& = & e_{2}  \\
\end{array}
&
\begin{array}{lll}
m(e_{2}, e_{2}, e_{1})& = & e_{1}+e_{2} \\
m(e_{2}, e_{2}, e_{2})& = & e_{1}+2e_{2} \\
m(e_{1}, e_{2}, e_{1})& = & e_{2} \\
m(e_{2}, e_{1}, e_{2})& = & e_{1}+e_{2}
\end{array}%
\end{array}
\end{equation*}

The following matrices correspond to automorphisms of the previous totally associative ternary algebra.

\begin{center}
\begin{tabular}{lll}
  $\left(\begin{array}{cc}
1 & 0 \\
0 & 1%
\end{array} \right)$&$ \left(
\begin{array}{cc}
-1 & 0 \\
0 & -1%
\end{array}%
\right)$ & $\ \ \left(
\begin{array}{cc}
-1 & -1 \\
0 & 1%
\end{array}%
\right) $ \\
  $\left(
\begin{array}{cc}
-\frac{1}{\sqrt{5}} & -\frac{3}{\sqrt{5}} \\
\frac{2}{\sqrt{5}} & \frac{1}{\sqrt{5}}%
\end{array}%
\right) $ & $\  \left(
\begin{array}{cc}
1 & 1 \\
0 & -1%
\end{array}%
\right) $ &$\left(
\begin{array}{cc}
\frac{1}{\sqrt{5}} & \frac{3}{\sqrt{5}} \\
-\frac{2}{\sqrt{5}} & -\frac{1}{\sqrt{5}}%
\end{array}%
\right) $ \\
\end{tabular}
\end{center}

We deduce for example using Theorem \ref{thm2} the following
2-dimensional totally Hom-associative ternary algebras
$(\K^2,\widetilde{m}_i,\rho_i )$ for $i=1,2$, which are not totally
associative. They are  given with respect to a basis $\{ e_1,e_2\}$, by the
following  product
\begin{equation*}
\begin{array}{ll}
\begin{array}{lll}
\widetilde{m}_1(e_{1}, e_{1}, e_{1})& = & e_{1}\\
\widetilde{m}_1(e_{1}, e_{1}, e_{2})& = & e_1- e_2 \\
\widetilde{m}_1(e_{1}, e_{2}, e_{2})& = & 2 e_1- e_2 \\
\widetilde{m}_1(e_{2}, e_{1}, e_{1})& = & e_1- e_2  \\
\end{array}
&
\begin{array}{lll}
\widetilde{m}_1(e_{2}, e_{2}, e_{1})& = & 2 e_1- e_2 \\
\widetilde{m}_1(e_{2}, e_{2}, e_{2})& = & 3e_{1}-2 e_2 \\
\widetilde{m}_1(e_{1}, e_{2}, e_{1})& = & e_1- e_2 \\
\widetilde{m}_1(e_{2}, e_{1}, e_{2})& = & 2e_1- e_2
\end{array}%
\end{array}
\end{equation*}

and the linear map
\begin{eqnarray*}
\rho_1 (e_1)&=& e_1 \\
\rho_1 (e_2)&=& e_1- e_2
\end{eqnarray*}

The second example is given by

\begin{equation*}
\begin{array}{ll}
\begin{array}{lll}
\widetilde{m}_2(e_{1}, e_{1}, e_{1})& = &
\frac{1}{\sqrt{5}} e_1- \frac{2}{\sqrt{5}}e_2\\
\widetilde{m}_2(e_{1}, e_{1}, e_{2})& = &
\frac{3}{\sqrt{5}} e_1- \frac{1}{\sqrt{5}}e_2 \\
\widetilde{m}_2(e_{1}, e_{2}, e_{2})& = & \frac{4}{\sqrt{5}} e_1-
\frac{3}{\sqrt{5}}e_2\\
\widetilde{m}_2(e_{2}, e_{1}, e_{1})& = &
\frac{3}{\sqrt{5}} e_1- \frac{1}{\sqrt{5}}e_2  \\
\end{array}
&
\begin{array}{lll}
\widetilde{m}_2(e_{2}, e_{2}, e_{1})& = &
\frac{4}{\sqrt{5}} e_1- \frac{3}{\sqrt{5}}e_2 \\
\widetilde{m}_2(e_{2}, e_{2}, e_{2})& = & \frac{7}{\sqrt{5}}
e_1-\frac{4}{\sqrt{5}}e_2\\
\widetilde{m}_2(e_{1}, e_{2}, e_{1})& = &
\frac{3}{\sqrt{5}} e_1- \frac{1}{\sqrt{5}}e_2 \\
\widetilde{m}_2(e_{2}, e_{1}, e_{2})& = & \frac{4}{\sqrt{5}} e_1-
\frac{3}{\sqrt{5}}e_2
\end{array}%
\end{array}
\end{equation*}

and the linear map
\begin{eqnarray*}
\rho_2 (e_1)&=& \frac{1}{\sqrt{5}} e_1- \frac{2}{\sqrt{5}}e_2 \\
\rho_2 (e_2)&=& \frac{3}{\sqrt{5}} e_1- \frac{1}{\sqrt{5}}e_2
\end{eqnarray*}

\end{example}

\section{Tensor products } \label{sec3:tensorproducts}
In this section we define the tensor product of two $n$-ary
Hom-algebras and prove some results involving $n$-ary Hom-algebras
of Lie type and associative type.

Let $\mathcal{A}=(V,m,\alpha)$, where
$\alpha=(\alpha_{i})_{i=1,\cdots,n-1}$,  and
$\mathcal{A}'=(V',m',\alpha')$,
$\alpha'=(\alpha'_{i})_{i=1,\cdots,n-1}$, be two $n$-ary
Hom-algebras of given type. The tensor product $\mathcal{A}\otimes
\mathcal{A}'$ is an $n$-ary algebra defined by the triple $(V\otimes
V',m\otimes m',\alpha\otimes \alpha')$ where
$\alpha\otimes\alpha'=(\alpha_i\otimes\alpha'_{i})_{i=1,\cdots,n-1}$
with
\begin{eqnarray*}
m\otimes m'(x_{1}\otimes x'_{1},\cdots,x_{n}\otimes
x'_{n})&=&m(x_{1},\cdots,x_{n})\otimes
m^{\prime}(x'_{1},\cdots,x'_{n}),\\
\alpha_i\otimes\alpha'_{i}(x\otimes
x')&=&\alpha_i(x)\otimes\alpha'_{i}(x'),
\end{eqnarray*}
where $x,x_{1},\cdots,x_{n}\in V$ and   $x',x'_{1},\cdots,x'_{n}\in
V'$.

Recall that an $n$-ary multiplication $m$ is symmetric if
\begin{equation}
m(x_{\sigma (1)}, \cdots, x_{\sigma (n)})=m(x_{1}, \cdots,
x_{n}),\quad \forall \ \sigma \in \mathcal{S}_{n}\text{ and }
\forall \ x_{1},\cdots,x_{n}\in V
\end{equation}%
and $m'$ is skew-symmetric if
\begin{equation}
m'(x_{\sigma (1)}, \cdots, x_{\sigma (n)})=Sgn(\sigma )m'(x_{1},
\cdots, x_{n}),\quad \forall \ \sigma \in \mathcal{S}_{n}
\text{ and }\forall \ x_{1},\cdots,x_{n}\in V'
\end{equation}%
where $\mathcal{S}_{n}$ stands for the
permutation group on $n$ elements.

We have the following obvious lemma.
\begin{lem}
Let $\mathcal{A}=(V,m,\alpha)$  and $\mathcal{A}'=(V',m',\alpha')$
be two $n$-ary Hom-algebras of given type. If $m$ is symmetric and
$m'$ is skew-symmetric then $m\otimes m'$ is skew-symmetric.
\end{lem}
\begin{proposition}
Let $\mathcal{A}$ be an $n$-ary totally Hom-associative algebra and
$\mathcal{A}'$ be an $n$-ary partially Hom-associative algebra, then
 $\mathcal{A}\otimes \mathcal{A}'$ is an $n$-ary partially Hom-associative algebra.
\end{proposition}

\begin{proof}
 Using the   $n$-ary totally Hom-associativity of $\A$ and  $n$-ary partially Hom-associativity of $\A '$, we obtain for $x_{1},\cdots,x_{2n-1}\in \mathcal{A}$ and $x'_{1},\cdots,x'_{2n-1}\in \mathcal{A}'$.
\begin{eqnarray*}
\sum_{i=0}^{n-1}{
m\otimes m'(\alpha_1\otimes\alpha'_{1}(x_1\otimes
x'_1),\cdots,\alpha_i\otimes\alpha'_{i}(x_i\otimes
x'_i),m\otimes m'(x_{i+1}\otimes x'_{i+1},\cdots,x_{i+n}\otimes x'_{i+n}),}\\
\alpha_{i+1}\otimes\alpha'_{i+1}(x_{i+1+n}\otimes x'_{i+1+n}),\cdots ,\alpha_{n-1}\otimes\alpha'_{n-1}(x_{2n-1}\otimes x'_{2n-1}))\\
=\sum_{i=0}^{n-1}{
m(\alpha_1(x_1),\cdots,\alpha_i(x_i),m(x_{i+1}\cdots,x_{i+n}),\alpha_{i+1}(x_{i+1+n})
,\cdots ,\alpha_{n-1}(x_{2n-1}))}
\otimes \\ m'(\alpha'_{1}(x'_1),\cdots,\alpha'_{i}(x'_i), m'
(x'_{i+1},\cdots ,x'_{i+n}),
\alpha'_{i+1}( x'_{i+1+n}),\cdots ,\alpha'_{n-1}( x'_{2n-1}))\\
=
m(\alpha_1(x_1),\cdots,\alpha_j(x_j),m(x_{j+1}\cdots,x_{j+n}),\alpha_{j+1}(x_{j+1+n})
,\cdots ,\alpha_{n-1}(x_{2n-1}))
\otimes \\\sum_{i=0}^{n-1}{ m'(\alpha'_{1}(x'_1),\cdots,\alpha'_{i}(x'_i), m'
(x'_{i+1},\cdots ,x'_{i+n}),
\alpha'_{i+1}( x'_{i+1+n}),\cdots ,\alpha'_{n-1}( x'_{2n-1}))}=0
\end{eqnarray*}
\end{proof}

\begin{proposition}
Let $\mathcal{A}$ be a symmetric $n$-ary totally Hom-associative
algebra and $\mathcal{B}$ be an $n$-ary Lie algebra, then
 $\mathcal{A}\otimes \mathcal{B}$ is an $n$-ary Lie algebra.
\end{proposition}

\begin{proof}
The proof is similar to the proof of the previous proposition.
\end{proof}

\end{document}